\documentclass[12pt]{article}
\usepackage{amsmath,amsfonts}

\usepackage{palatino}
\usepackage{graphicx}
\usepackage{epstopdf}
\usepackage[sort,compress]{cite}

\newcommand{\eps}{{\varepsilon}}

\newcommand{\R}{\mathbb{R}}

\title{Small-noise analysis and symmetrization of implicit Monte Carlo
  samplers}

\date{October 21, 2014}

\author{Jonathan Goodman\thanks{Department of Mathematics, Courant
    Institute of Mathematical Sciences, New York University, New York,
    NY 10012, USA}\\
  Kevin K.~Lin\thanks{School of Mathematics, University of Arizona,
    Tucson, AZ 85721, USA}\\
  Matthias Morzfeld\thanks{Department of Mathematics, University of
    California at Berkeley and Lawrence Berkeley National Laboratory,
    Berkeley, CA 94720, USA}}

\begin{document}

\maketitle

\begin{abstract}
\noindent
Implicit samplers are algorithms for producing independent, weighted
samples from multi-variate probability distributions.  These are often
applied in Bayesian data assimilation algorithms.  We use Laplace
asymptotic expansions to analyze two implicit samplers in the small
noise regime.  Our analysis suggests a symmetrization of the algorithms
that leads to improved (implicit) sampling schemes at a relatively small
additional cost.  Computational experiments confirm the theory and show
that symmetrization is effective for small noise sampling problems.
\end{abstract}

\section{Introduction} \label{sec:intro}

Markov chain Monte Carlo (MCMC) techniques are widely
used for sampling complicated distributions.  However, some data
assimilation methods rely on independent samples from known
distributions \cite{Bocquet2010,PeterJan2009,Fournier2010,Doucet2001}.
Weighted direct samplers give independent samples from a proposal
distribution that is not the target distribution, and compensate for
this with a random weight factor (see e.g.~\cite{ChorinHald} and
references there).  The variance of the weight factor determines the
quality of the sampler \cite{liuchen1995,GordonReview,Doucet2001}.

This paper studies two weighted direct samplers.

One, the {\em linear map} method, has been proposed independently
several times and has several names in the literature; see, e.g.,
\cite{Chorin2010}, and also \cite{An} for a similar method.  The other
is the {\em random map} method, which was proposed in
\cite{Morzfeld2011}.  The linear and random map methods can both be
viewed as examples of \emph{implicit samplers}
\cite{chorintupnas,atkins,Chorin2010}.  We introduce a small noise
parameter, $\eps$, similar to that of \cite{Weare2012,Weare2013}, and
analyze the performance of these algorithms in sampling general smooth
probability densities on finite-dimensional spaces in the limit $\eps
\to 0$.  The methods we study use a Gaussian approximation to the target
distribution, which is valid in the small noise limit.  Many data
assimilation applications are in the small noise regime.

We study a standard quality measure of weighted direct samplers.  Our
analysis consists in calculating {\em error constants}, which are the
coefficients of the leading powers of $\eps$ in the small noise
asymptotic expansion of the quality measures.  As long as simple
smoothness hypotheses are satisfied, the error constants for the linear
and random map methods differ by a factor that depends only on the
dimension.  This factor converges to one as the dimension goes to
infinity.

The form of the error constant suggests that a symmetrization may remove
the leading error term.  We study symmetrized versions of the linear and
random map methods to confirm this.  The error is one order smaller in
$\eps$.  The error constants are not exactly proportional, but their
ratio does converge to one as the dimension converges to infinity.  We
present computational experiments that confirm the small noise
asymptotic calculations.  The numerical experiments further demonstrate
that the symmetrized methods are more accurate in the small noise
regime, and show that the symmetrized methods may perform significantly
better than the corresponding ``simple'' methods even when the noise is
not so small.

\medskip

This paper is organized as follows: in Section~\ref{sect:main-results},
we set up the notation, present the algorithms, explain how they can be
symmetrized, and summarize our theoretical results.
Section~\ref{sect:preliminaries} describes two general technical tricks
that simplify the asymptotic analysis.  The relatively simple
derivations of the linear map results are in Section
\ref{sec:linear-map}.  These depend on Subsection
\ref{sec:variance-lemma} only.  The ideas in Subsection \ref{sec:Wick}
are needed only for the explicit error constant formulas for the first
computational example in Section \ref{sec:computations}.  The analysis
in Section \ref{sec:random-map} of random map methods also uses the
formula derived in Subsection \ref{sec:trick}.  Section
\ref{sec:computations} describes numerical experiments on two test
problems which confirm the asymptotic theory in detail.  It may be read
without the theoretical Sections \ref{sect:preliminaries},
\ref{sec:linear-map}, and \ref{sec:random-map}.
Section~\ref{sec:Conclusions} summarizes our views of these results and
puts them in context.

\section{Algorithms, symmetrization and main results}
\label{sect:main-results}

In this section, we describe the sampling algorithms to be studied in
the rest of the paper.  We introduce the small-noise scaling used in the
analysis in Section~\ref{sect:Summary of main results}, where we also
state our main theoretical scaling results.

The following notation is used throughout the paper.  Let $f(x)$ and
$r(x)$ be two functions on $\R^d~.$ We write $f \propto r$ if $f(x) = C
r(x)$ for some fixed $C$.  If distributions depend on $\eps$, we write
$f(x,\eps) \propto r(x,\eps)$ if there is a $C_{\eps}$ with $f(x,\eps) =
C_{\eps}r(x,\eps)$.  For any non-negative $f$ with $0<\int
f(x)~dx<\infty~,$ there is a probability density $p \propto f$.  We
write $X \sim p$ if $X$ is a random variable whose probability density 
is $p$.  We say a random variable $X \sim q$ together with a non-negative
weight function $w$ is a {\em weighted sample} of a given probability
density $p$ if
\begin{equation}
  E_p(u(X)) = \frac{E_q(u(X)\cdot w(X))}{E_q(w(X))},
\label{eqn:weighted-sample}  
\end{equation}
for every bounded continuous function $u$.  A {\em weighted sampler} of
$p$ with {\em proposal} $q$ is a stochastic algorithm that produces $X
\sim q \propto g$.  It is a {\em direct sampler} if successive samples
are independent.  The direct samplers we consider have deterministic
weight functions
\begin{equation}
   w(x) = \frac{f(x)}{g(x)}\propto \frac{p(x)}{q(x)}  ~.
\label{eqn:w} \end{equation}
We assume $f$ and $g$ may be evaluated, but the normalizing
constants may be unknown (as is typically the case in applications).

A perfect sampler would have a constant weight function $w =C$, which
would force $p=q$.  We measure the quality of a weighted sampler by the
non-dimensionalized deviation of $w$ from a constant:
\begin{equation}
R := \frac{E_q(w(X)^2)}{E_q(w(X))^2}~,~~~ Q := R-1 \; .
\label{eqn:Q} \end{equation}
The quality measure $Q$ was also used in \cite{Weare2013},
and several motivations for it are given in 
\cite{Doucet2001,GordonReview,LiuChen98,Bergmman99}.  
In particular, a heuristic relates a collection of $N$ independent weighted
samples to $N/R$ independent un-weighted samples, making $N/R$ an
effective sample size. 
A small $Q$ is important 
in recursive particle filter algorithms.
There, probability densities are sampled recursively,
as the data are collected,
and the weights accumulate as a product of 
the weights at each step.
Thus, the $R$ of the product can grow rapidly 
if the $Q$ of each of the factors is not small.
Both algorithms analyzed here have the property that $Q\to0$ as
$\eps\to0~;$ the question that concerns us is the rate of convergence.

The methods we consider sample 
\begin{equation}
  p(x) \propto f(x) = e^{-F(x)}~.
  \label{eqn:p-no-eps}  \end{equation}
We assume that $F$ is smooth and has a single global minimum, which is
non-degenerate.  Unless otherwise stated,
we also assume (without loss
of generality) that the minimum is located at $x=0~,$ and that $F(0)=0$.
We write a Taylor expansion of $F$ near zero as
\begin{equation}
F(x) = {\textstyle \frac{1}{2}}x^tHx + C_3(x) + \cdots + C_6(x) + O(\left|x\right|^{7}) ~,
\label{eqn:FTs}  \end{equation}
where $H$ is the Hessian matrix of $F$ at $x=0$, 
and $C_k(x)$ is the homogeneous polynomial of 
degree $k$ 
\begin{equation}
C_k(x) =  \frac{1}{k!} \sum_{|\alpha|=k} x^{\alpha} \partial_x^{\alpha} F(0)  ~.
\end{equation}
The random map samplers also require that certain equations related to $F$
have unique and well behaved solutions (see below).

\subsection{Simple and symmetrized linear map methods}

The simple linear map method uses $X \sim \pi$, where $\pi$ is the local 
Gaussian approximation that uses the first term on the right of (\ref{eqn:FTs}),
\begin{equation}
\pi(x) \propto e^{-x^tHx/2} ~.
\label{eqn:pie} \end{equation}
Direct Gaussian sampling algorithms make this possible.
Using (\ref{eqn:w}) with $f=e^{-F}$ and $g=e^{-x^tHx/2}$,
we find the weight function
\begin{equation}
  w(x) = e^{-F(x) + x^t H x/2}~.
  \label{eqn:LinearMapWeight}
\end{equation}
The simple linear map Monte Carlo algorithm to estimate $E_p(u(X))$ is: 
\begin{enumerate}{\roman{enumi}}
\item Generate $N$ independent Gaussian samples $X_k \sim \pi$
\item Compute weights $W_k = w(X_k)$ using (\ref{eqn:LinearMapWeight})
\item Compute the estimator
  \begin{displaymath}
    \frac{\sum_{k=1}^N u(X_k)w(X_k)}{\sum_{k=1}^N w(X_k)},
  \end{displaymath}
  of $E_p(u(X))~.$
\end{enumerate}
.

In practice, the minimizer of $F$ will not be at $0~,$ 
nor will its Hessian at the minimum be the identity.
It will be necessary to first find $x_* =
\mbox{argmin }F(x)$ and evaluate $H(x_*)$, 
the Hessian of $F$ at the minimum. 
This can be a time-consuming step.

As we will see below, the leading-order term in the $\eps$-expansion of
$Q$ depends only on $C_3$ (see equation (\ref{eqn:A-linear-map})).  This
is not surprising, as the simple method is based on the approximation
$F(x) \approx \frac{1}{2}x^tHx$, and $C_3(x)$ is the largest correction.
Since $C_3$ is an odd function of $x$, one may hope that the leading
error term can be removed by a symmetrization related to the classical
Monte Carlo trick of antithetic variates \cite{HammersleyHandscomb,
KalosWhitlock}.  Here we present a symmetrized linear map method; we
will verify in Section~\ref{sec:linear-map-symm} that it removes the
principal error term in the small noise limit.

The symmetrized linear map method is as follows: first draw $\xi
\sim \pi$ as before, and evaluate the linear map weights 
(\ref{eqn:LinearMapWeight}) for $\xi$ and for
$-\xi$.  Note that $\pi(\xi) = \pi(-\xi)$ so these weights are
\begin{equation}
  w_+ = \frac{f(\xi)}{\pi(\xi)} \;\; , \;\;\;\;
  w_- = \frac{f(-\xi)}{\pi(\xi)}.
\end{equation}
Then return $X=\xi$ or
  $X = -\xi$ with probabilities
\begin{equation}
p_+ = \frac{w_+}{w_++w_-} \;\; , \;\;\;\; p_- = \frac{w_-}{w_++w_-} ~.
\label{eqn:wpm}  \end{equation}
These probabilities have a particle filter interpretation.
Consider $(\xi,w_+)$ and $(-\xi,w_-)$ to form a two element weighted ensemble.
The formulas (\ref{eqn:wpm}) are the probabilities that would be used to sub-sample this to 
a one element un-weighted ensemble \cite{GordonReview}.

To find the weight function of the symmetrized linear map method we must
identify $q_s(x)$, the probability density of $X$.  There are two ways
to generate $X=x$ (using the convention that $x$ is a possible value
of the random variable $X$). One way is to propose
$\xi=x$ and then take the $+$ choice in (\ref{eqn:wpm}).  The other way
is to propose $\xi=-x$ and then take the $-$ choice.  The probability
density for $\xi=x$ is $\pi(x)$.  The density for $-\xi$ is $\pi(-\xi) =
\pi(\xi)$.  The probability to get $x$ if $x$ was proposed is
\begin{equation}
p_+(x) = \frac{w(x)}{w(x) + w(-x)} \; .
\end{equation}
The probability to get $x$ if $-x$ was proposed is the same, since
\begin{align}
  p_-(-x) &= \frac{w_-(-x)}{w_+(-x) + w_-(-x)} ,\nonumber\\[1ex]
  &= \frac{w(-(-x))}{w(-x) + w(-(-x))} ,\nonumber\\[1ex]
  &= p_+(x) \; .
\end{align}
Therefore, the pdf of $X$ is
\begin{align}
  q_s(x) &= \pi(x) p_+(x) + \pi(-x) p_-(-x),\nonumber\\[1ex]
  &= \pi(x) \frac{2w(x)}{w(x) + w(-x)} ~.
\label{eqn:qs}  
\end{align}
Moreover, if $\pi(x)$ is a normalized probability density, 
then $q_s$ is also normalized.
This can be seen by using the right side  of (\ref{eqn:qs})
\begin{align}
  \int q_s(x) \,dx &= {\textstyle \frac{1}{2}}\int \left( q_s(x) + q_s(-x) \right) \,dx ,\nonumber\\
  &= {\textstyle \frac{1}{2}}\int \frac{2\pi(x)}{w(x) + w(-x)} \left(w(x) +w(-x) \right) \,dx, \nonumber\\
  &=1 ~.
\end{align}
The weight function (\ref{eqn:w}) for the symmetrized method is thus
\begin{equation}
  w_s(x) \propto \frac{p(x)}{q_s(x)} \propto \frac{\pi(x) w(x)}
                                    {\pi(x) \frac{ 2w(x)}{w(x) + w(-x)}}
         = \frac{w(x) + w(-x)}{2}  ~.
\label{eqn:ws}  
\end{equation}

The simple linear map sampler and the symmetrized sampler have
different symmetries.  The simple sampler has a symmetric proposal
density and non-symmetric weight.  The symmetrized sampler has a
non-symmetric proposal density, $q_s(x) \neq q_s(-x)$, but a symmetric
weight function.
Intuitively one can therefore expect that the quality measure of the
symmetrized method is better because a symmetric weight function is
``more nearly constant'' for small $x$, particularly in the small noise 
regime described below.

\subsection{Simple and symmetrized random map methods}
\label{sect:Random map methods}

The {\em simple} (as opposed to {\em symmetrized}) random map method is
described in \cite{Morzfeld2011}.  We review the method here for
notation and completeness.

One first samples $\xi \sim \pi$, and then chooses 
\begin{equation}
  X = \lambda(\xi) \xi~.
  \label{eqn:random-map-stretch} 
  \end{equation}
The stretch factor $\lambda(\xi)\geq 0$ is defined implicitly via
\begin{equation}
F(\lambda(\xi)\xi) = {\textstyle\frac{1}{2}} \xi^tH\xi ~.
\label{eqn:stretch-condition} \end{equation}
The random map algorithm gets its name from the map $\xi \mapsto X~.$ To
ensure the correctness of the algorithm, 
we need to assume that equation~(\ref{eqn:stretch-condition}) 
has a unique solution $\lambda>0$ for every $\xi\neq0~.$ 
This will be the case if, e.g., every level set (except the zero level set) 
of $F$ is ``star-shaped,'' i.e.,
for every $c>0~,$ every straight line through 0 intersects the level set
$F^{-1}(c)$ transversely at exactly two points.

To determine the weight function of the random map method,
note that if $\xi \sim \pi$ and
$X = x(\xi)$, then $X$ has probability density 
\begin{equation}
   q(x(\xi)) = \pi(\xi) \left|\det\!\left(\frac{\partial \xi}{\partial x}\right)\right|,
\end{equation}
so that we find the weight $w$ from (\ref{eqn:w}) to be
\begin{equation}
  w(\xi) = \left|\det\!\left(\frac{\partial x}{\partial \xi}\right)\right| ~,
\end{equation}
choosing the arbitrary implicit constant to be equal to 1 here.

The Jacobian determinant is
\begin{equation}
   w(\xi) 
      = \lambda(\xi)^{d-1} \frac{\xi^tH\xi}{\left|\xi^t \nabla_x F(\lambda(\xi)\xi)\right|} ~.
\label{eqn:w-random-map}  \end{equation}
To see this, note that the Jacobian matrix is obtained by  
differentiating (\ref{eqn:random-map-stretch})
\begin{equation}
  \frac{\partial x}{\partial \xi} =\xi \, [\nabla \lambda(\xi)]^t + \lambda(\xi)I,
\end{equation}
where $\nabla \lambda$ is the column vector with entries 
  $\partial_{\xi_j} \lambda(\xi)$.
The determinant identity 
\begin{equation}
  \mbox{det}(\lambda I + A ) = \lambda^d + \lambda^{d-1} \mbox{tr}(A) + \cdots \;,
\end{equation}
gives
\begin{equation}
  \mbox{ det}\!\left(\frac{\partial x}{\partial \xi}\right) 
   = \lambda^d + \lambda^{d-1} \xi^t \nabla \lambda ~,
\end{equation}
where the terms of order $\lambda^{d-2}$ 
and lower vanish because $ \xi\,  [\nabla \lambda(\xi)]^t$ 
is a matrix of rank one.
A calculation (given just below) gives
\begin{equation}
  \xi^t \nabla_{\xi} \lambda = \lambda \left( \frac{\xi^tH\xi}{x^t \nabla_x F} -1 \right),
\label{eqn:xi-dot-grad-lambda} \end{equation}
which immediately leads to (\ref{eqn:w-random-map}).
To verify (\ref{eqn:xi-dot-grad-lambda}), 
we differentiate (\ref{eqn:stretch-condition}) 
with respect to $\xi_i$:
\begin{align}
  \sum_j \partial_{x_j} F(\lambda(\xi)\xi)\, \partial_{\xi_i} [\lambda(\xi) \xi_j] 
          &= (H\xi)_i \;, \\
  \sum_j \partial_{x_j} F(\lambda(\xi)\xi)\, 
         \left( \frac{\partial \lambda(\xi)}{\partial \xi_i}\xi_j + \lambda(\xi) \delta_{ij} \right)
          &= (H\xi)_i \; .
\end{align}
We multiply by $\xi_i$, sum over $i$, and use the relations $\lambda(\xi)\xi = x$, and 
$\xi = \frac{1}{\lambda(\xi)}x$:
\begin{align}
     \xi^t \nabla_x F(\lambda(\xi)\xi) \,\xi^t \nabla_{\xi} \lambda(\xi) 
   + \lambda(\xi) \xi^t \nabla_x F(\lambda(\xi)\xi) 
 &= \xi^tH\xi, \\
     \frac{1}{\lambda(\xi)} x^t \nabla_x F(\lambda(\xi)\xi) \, \xi^t \nabla_{\xi} \lambda(\xi)
   + x^t \nabla_x F(\lambda(\xi)\xi) 
 &=  \xi^tH\xi \;.
\end{align}
Solving for $\xi^t\nabla_x \lambda(\xi)$ gives (\ref{eqn:xi-dot-grad-lambda}).

Our symmetrization of the simple random map method is a natural
adaptation of the symmetrization of the simple linear map method.
There are three steps:
\begin{enumerate}{\roman{enumi}}
\item Generate a sample $\xi\sim \pi$.
\item Compute $x_+ = \lambda(\xi)\cdot\xi$ and $x_- =
  \lambda(-\xi)\cdot(-\xi)$, each using (\ref{eqn:stretch-condition}).
\item Use $x = x_+$ with probability $p_+(\xi) := w(\xi)/(w(\xi)+w(-\xi))$.
  Otherwise use $x_-$.
\end{enumerate}
The arguments leading to (\ref{eqn:qs}) and (\ref{eqn:ws}) apply here too.
The probability density of $X$ produced by the symmetrized random map
method is thus
\begin{equation}
q_s(x) = \frac{2}{w(\xi(x))+w(-\xi(x))} e^{-F(x)}~,
\label{eqn:symmetrized-random-map-proposal-pdf}
\end{equation}
where $w(\xi)$ is the weight of the simple random map method.
The weight function for the symmetrized random map method is
\begin{equation}
  w_s(\xi) = \frac{w(\xi)+w(-\xi)}2~,
  \label{eqn:symmetrized-random-map-weights}
\end{equation}
where $w$ is the weight of the simple method (\ref{eqn:w-random-map}).

\subsection{Summary of small noise theory}
\label{sect:Summary of main results}

The small noise problem concerns the scaled density
\begin{equation}
  p(x) \propto f(x) = e^{-F(x)/\eps}~.
  \label{eqn:p}  \end{equation}
Recursive particle filter applications often call for proposal distributions 
roughly of the form (\ref{eqn:p}). 
When the noise parameter $\eps$ is small, most of the probability in $p$ is near the 
point of maximum probability, which we continue to take to be $x_*=0$.
Therefore $F(x) \approx \frac{1}{2}x^tHx$ (see (\ref{eqn:FTs})) may be a useful approximation.

We state and derive the small noise theory using a standard scaling,
\begin{equation}
	\widetilde{x} = \eps^{1/2} x ~.
\label{eqn:x-tilde} 
\end{equation}
This scales the terms in the Taylor expansion (\ref{eqn:FTs}) as
\begin{equation}
   \frac{F(\widetilde{x})}{\eps} = {\textstyle \frac{1}{2}} \widetilde{x}^tH\widetilde{x}
      + \eps^{1/2}C_3(\widetilde{x}) + \eps C_4(\widetilde{x}) 
      + \eps^{3/2}C_5(\widetilde{x}) + \eps^{2}C_6(\widetilde{x}) + O(\eps^{5/2}) ~.
\label{eqn:F-scaled}  \end{equation}
The target density therefore satisfies
\begin{equation}
  p(\widetilde{x}) \propto \exp\left(-\widetilde{x}^tH\widetilde{x}/2-\eps^{1/2}C_3(\widetilde{x}) 
                 - \eps C_4(\widetilde{x}) 
        - O(\eps^{3/2})\right).
\label{eqn:p-pie}  
\end{equation}
For the rest of the theory, we assume $p$ satisfies (\ref{eqn:p-pie}).
Following common practice, we drop the tilde.

Our theoretical results take the form of asymptotic approximations of $Q$ defined in (\ref{eqn:Q}).
The simple linear map and random map methods have the scaling  
\begin{equation}
Q = \eps A + O(\eps^{3/2}).
\label{eqn:root-eps} \end{equation}
The error constants are
\begin{align}
A &=& &E_{\pi}(C_3(X)^2) &\qquad\mbox{(simple linear map),} \label{eqn:A-linear-map}\\
A &=& \frac{(1+d)^2}{(2+d)(4+d)}&E_{\pi}(C_3(X)^2) &\qquad\mbox{(simple random map).} 
           \label{eqn:A-random-map}
\end{align}
The error scaling for the symmetrized methods is
\begin{equation}
Q =  \eps^2 B + O(\eps^{5/2}) ~,
\label{eqn:eps-symm} \end{equation}
with error constants of the form
\begin{align}
B &=   \mbox{var}_{\pi}\!\left(\, C_4 - {\textstyle \frac{1}{2}}C_3^2 \right)   &\mbox{(symmetrized linear map),} 
                                                                 \label{eqn:A-linear-map-symm}\\
B &=   \mbox{var}_{\pi}\!\left(\, C_4 - {\textstyle \frac{1}{2}}C_3^2 \right) +
    c_d\cdot K  &\mbox{(symmetrized random map).}  
                                                                 \label{eqn:A-random-map-symm} 
\end{align}
Here  $c_d = O(1/d)$, and $K$ is a
possibly dimension-dependent constant depending on $F$.
The exact form is given in Section~\ref{sec:symmetrized-random-map}.

We have the following conclusions.  For both methods, $Q\to 0$ in the
small noise limit $\eps\to 0$.  This is perhaps not surprising because
the Gaussian approximation $\pi$ becomes exact in this limit.  On the
other hand, this property cannot be taken for granted in general; see,
e.g., \cite{Weare2012}.

The simple linear and random map methods 
have the same order as $\eps \to 0$,
and the symmetrized methods have a higher order.
Thus, for any fixed problem and for sufficiently small $\eps$, 
the error constants of the symmetrized methods are significantly smaller
than the error constants for the corresponding simpler methods.
The ratio of the error constants for the linear and random map methods
depends only on the dimension.
This factor converges to 1 as $d \to \infty$.
Thus, in the limits $\eps \to 0$ and $d \to \infty$, 
the random map methods lose their advantages over the linear map
methods.

\section{Analysis tools}
\label{sect:preliminaries}

Here we describe two tools that we will use in the analysis
of the linear and random map methods.

\subsection{The variance lemma} \label{sec:variance-lemma}
The variance lemma is a simple way to understand some cancellations that
occur in computing $Q$ for small $\eps$.
It applies to functions $u(x,\eps)$ of the form
\[
u(x,\eps) = 1 + \eps^r u_1(x) + \eps^{2r} u_2(x) + O(\eps^{3r}) ~.
\]
It states that if 
\begin{equation}
   Q = \frac{ E(u(X,\eps)^2)}{E(u(X,\eps))^2} -1 ~,
\label{eqn:Q-general} \end{equation}
then
\begin{equation}
   Q = \eps^{2r}\mbox{var}(u_1(X)) + O(\eps^{3r}) ~.
\label{eqn:Q-var}  \end{equation}
The variance formula (\ref{eqn:Q-var}) does not depend on the distribution of $X$, 
except that the same distribution must be used throughout.
Expectations of $u_2$ appear at $O(\eps^{2r})$ in the numerator and denominator of 
(\ref{eqn:Q-general}) but they cancel in the ratio to leading order.

The verification is straightforward.
The numerator in (\ref{eqn:Q-var}) is
\begin{align*}
E(u(X,\eps)^2)  &= E( 1 + 2\eps^r u_1 + \eps^{2r}\left( u_1^2 + 2u_2\right) + O(\eps^{3r}), \\
            &= 1 + 2\eps^r E(u_1)+ \eps^{2r}\left( E(u_1^2) + 2 E(u_2)\right) + O(\eps^{3r}) ~.
\end{align*}
The denominator is
\begin{align*}
E(u(X,\eps))^2 &= \left( 1 + \eps^r E(u_1) + \eps^{2r}E(u_2) + o(\eps^{2r})\right)^2, \\
  &= 1 + 2 \eps^r E(u_1) + \eps^{2r} \left( E(u_1)^2 + 2 E(u_2)\right) + O(\eps^{3r}) ~.
\end{align*}
Therefore,
\begin{align*}
Q &= \frac{1 + 2\eps^r E(u_1) 
                + \eps^{2r}\left( E(u_1^2) + 2E(u_2)\right) + O(\eps^{3r})}{1 + 2 \eps^r E(u_1) + \eps^{2r} \left( E(u_1)^2 + 2 E(u_2)\right) + O(\eps^{3r}) }
                -1 \\[1ex]
   &= \left[ 1 + 2\eps^r E(u_1) 
                + \eps^{2r}\left( E(u_1^2) + 2E(u_2)\right) + O(\eps^{3r})\right] \\[1ex]
    &\;\; \times  \left[  1 - 2 \eps^r E(u_1) + \eps^{2r} 
            \left( 3E(u_1)^2 - 2 E(u_2)\right) + O(\eps^{3r})\right] -1 \\
      &= \eps^{2r} \left( E(u_1^2) - E(u_1)^2 \right) + O(\eps^{3r}) ~.
\end{align*} 
Note that this conclusion depends on the existence of a function $u_2(x)$, but it does not
depend on what $u_2$ is.

\subsection{Evaluating rational Gaussian expectations}
\label{sec:trick}

The random map analysis in Section~\ref{sec:random-map}
leads to Gaussian expectations of the form
\[
    E_{\pi}\left( \frac{C(\xi)}{|\xi|^{2r}}\right) ~, 
\]
where $C$ is a homogeneous polynomial of some degree.
These are related to expectations of $C$.
In fact, if $f(\xi)$ is homogeneous of degree $q$, then
\begin{equation}
   (q+d)E_{\pi}\!\left( f(\xi)\right)  = 
    E_{\pi}\!\left( |\xi|^2 f(\xi) \right)~ .
\label{eqn:Gauss-integral-identity}  \end{equation}
Taking $C$ of degree $p$, and $f(\xi)=C(\xi)/|\xi|^2$ or $f = C(\xi)/|\xi|^4$, 
gives $q = p-2$ or $q=p-4$, and
\begin{equation}
    E_{\pi}\!\left( \frac{C(\xi)}{|\xi|^2}\right)  = 
    \frac{1}{p-2+d}E_{\pi}\!\left( C(\xi) \right)~ ,
\label{eqn:C2}  \end{equation}
or (iterating twice)
\begin{equation}
    E_{\pi}\!\left( \frac{C(\xi)}{|\xi|^4}\right)  = 
    \frac{1}{(p-4+d)(p-2+d)}E_{\pi}\!\left( C(\xi) \right)~ .
\label{eqn:C4}  \end{equation}
This result, which may be derived as a $\Gamma$ function identity, 
is surely not new.

We give an elementary derivation that uses
the function 
\[
  I(r) = \int f(r \xi) e^{-|\xi|^2/2}\,d\xi ~.
\]
On one hand, 
\[
  I(r) = r^q \int f(\xi) e^{-|\xi|^2/2}\,d\xi = r^q I(1)~.
\]
On the other hand, we can change variables with $r \xi = \eta$ to get
\[
  I(r) = \int f(\eta) e^{-|\eta|^2/(2r^2)} \frac{1}{r^d}\,d\eta ~.
\]
Now differentiate with respect to $r$ and set $r = 1$:
\begin{align*}
   q r^{q-1}I(1) &= I^{\prime}(r) \\
    &=  \frac{1}{r^3}\int f(\eta) |\eta|^2 e^{-|\eta|^2/(2r^2)} \,\frac{d\eta}{r^d}, \ \\
    &\;\;\;\;\;\; -\;
   d\int f(\eta) |\eta|^2 e^{-|\eta|^2/(2r^2)} \,\frac{d\eta}{r^{d+1}},  \\
q I(1) &= \int f(\xi) | \xi |^2 e^{-| \xi |^2/2} \,d \xi - d I(1) ,\\
(q+d) \int f(\xi) e^{-| \xi |^2/2}\,d \xi &= \int | \xi |^2 f(\xi) e^{-| \xi|^2/2}\,d \xi ~.
\end{align*}
This is the desired (\ref{eqn:Gauss-integral-identity}).

\section{Analysis of linear map methods} \label{sec:linear-map}

This section contains the calculations behind the results
(\ref{eqn:A-linear-map}) and (\ref{eqn:A-linear-map-symm}).  
We estimate the expectations required for $Q$ (see (\ref{eqn:Q})) 
using the Laplace  asymptotic expansion method, see, e.g., \cite{Murray}.  
The calculations are easy to justify if $F$ has a unique global minimum
and $F\to\infty$ rapidly enough as $|x|\to\infty$.

\subsection{Laplace asymptotics, simple linear map}

We wish to calculate the expected value of the weight
and the expected value of the square of the weight for 
the linear map method in (\ref{eqn:LinearMapWeight}).
We use the Taylor expansion of $F$ in (\ref{eqn:F-scaled})
to obtain a Taylor expansion of the weight 
\begin{align}
  w(x) &= e^{-\eps^{1/2}C_3(x) - \eps C_4(x) + O(\eps^{3/2})} \nonumber \\
       &= 1 - \left[ \eps^{1/2} C_3(x) + \eps C_4(x) \right] 
            + {\textstyle \frac{1}{2}} \left[  \eps^{1/2} C_3(x) \right]^2 + O(\eps^{3/2}) \nonumber \\
       &= 1 - \eps^{1/2} C_3(x) + \eps \left[  {
                   \textstyle \frac{1}{2}} C_3(x)^2 - C_4(x) \right] + O(\eps^{3/2}) ~.
  \label{eqn:w-expansion}
\end{align}
Recall that $C_3(x)$ is an odd function of $x$ and $\pi(x)$ is symmetric.
Therefore $E_{\pi}(C_3)=0$, and the variance lemma (\ref{eqn:Q-general}) 
with $r=1/2$ gives
\[
   Q = \eps \,\mbox{var}_{\pi}(C_3) + O(\eps^{3/2}) = \eps E_{\pi}(C_3^2) + O(\eps^{3/2})~.
\]
This is the desired result (\ref{eqn:A-linear-map}).

\subsection{Laplace asymptotics, symmetrized linear map} \label{sec:linear-map-symm}
We obtain the Taylor expansion of the weight of the symmetrized linear map method
from (\ref{eqn:ws}) and from the expansion of the weight of the simple linear map
method in (\ref{eqn:w-expansion}).
We note that the term that is anti-symmetric in $x$, which is $C_3(-x) = - C_3(x)$,
cancels, so that
\begin{equation}
w_s(x) \approx  1 +  \eps \left[  {\textstyle \frac{1}{2}} C_3^2 -C_4(x)\right] ~.
\label{eqn:SymmetrizeWeight}
\end{equation}
To apply the variance lemma,
we first show that
\begin{equation}
  Q ~~=~~ \frac{E_{q_s}(w_s(X)^2)}{E_{q_s}(w_s(X))^2}-1 
    ~~=~~ \frac{E_{\pi}(w_s(\xi)^2)}{E_{\pi}(w_s(\xi))^2} -1~.
\end{equation}
This shows that we can average over $\xi$ instead of $X$ when computing
the quality measure $Q$.  To see why, note that (\ref{eqn:qs}) implies
that for any function $u$,
\begin{equation}
    E_{q_s}(u(X)) = E_{\pi}\left(\frac{2w(\xi)}{w(\xi) + w(-\xi)} \,u(\xi)\right) ~.
\label{eqn:qs-pie}   \end{equation}
Together with (\ref{eqn:ws}), this implies that 
\begin{align}
  E_{q_s}(w_s(X)) &= E_{\pi}\left(\frac{2w(\xi)}{w(\xi) + w(-\xi)} \,w_s(\xi)\right)\nonumber\\
  &= E_{\pi}\left(\frac{2w(\xi)}{w(\xi) + w(-\xi)} \, \frac{w(\xi) + w(-\xi)}{2} \right)\nonumber\\
  &= E_{\pi}(w(\xi)) \label{eqn:Eqw-epiw}
\end{align}
The last equality follows from the symmetry of $\pi$.
Similar algebra and symmetry reasoning leads to
\begin{align}
  E_{q_s}(w_s(X)^2) 
  &= E_{\pi}\left(\frac{2w(\xi)}{w(\xi) + w(-\xi)} \, \left[\frac{w(\xi) + w(-\xi)}{2}\right]^2 \right)\nonumber\\
  &= E_{\pi}\left(w(\xi) \,\frac{w(\xi) + w(-\xi)}{2} \right)\nonumber\\
  &= E_{\pi}\left(w_s(\xi)^2 \right) ~.   \label{eqn:Eqwsq-epiwsq}
\end{align}
Application of the variance lemma to the above 
expression, with expectations over $\xi$,
and using (\ref{eqn:SymmetrizeWeight}),
leads to the error term (\ref{eqn:A-linear-map-symm}).

\subsection{Evaluating $E(C_3^2)$ with Wick's formula} 
\label{sec:Wick}

There is a more explicit expression for $E_{\pi}(C_3^2)$ based on Wick's
formula \cite{Koopmans}.  Recall that the distribution of a mean zero
multivariate Gaussian is completely determined by its covariance matrix.
Therefore, the expected value of a higher order monomial is a function
of the covariances.  Wick's formula (\ref{eqn:Wick}) is this function.
Of course, the expected value of an odd order monomial is zero; Wick's
formula gives the even-order moments.

The general version of Wick's formula is as follows (see, e.g., \cite{Koopmans}).
Suppose $X = (X_1,\ldots,X_d) \in {\mathbb R}^d$ is a multivariate mean zero Gaussian with 
covariances $C_{ij} = E_{\pi}(X_iX_j)$.
Let $i_k$, for $k = 1, \ldots,2n$, be a list of indices, with repeats allowed.
Let $M = X_{i_1} \cdots X_{i_{2n}}$ be the corresponding degree $2n$ monomial.
A  {\em pairing} is a partition of $\left\{1,\ldots,2n\right\}$ into $n$ sets of size $2$
\[
   P = \left\{ \left\{ k_i,l_1\right\}, \ldots, \left\{k_n,l_n\right\} \right\} ~.
\]
A pairing has the property that
\[
   \left\{ 1,\ldots,2n\right\} = \bigcup_{r=1}^n \left\{ k_r,l_r \right\} ~.
\]
The set of all pairings is $\cal P$.
The number of pairings is 
\[
    \left|{\cal P}\right| = (2n-1)(2n-3) \cdots 3 = (2n-1)!! ~.
\]
There are no pairings of a set with an odd number of elements.
Wick's formula gives the expected value of a monomial of even degree as a sum over
all pairings of the indices:
\begin{equation}
   E_{\pi}\left(  \prod_{k=1}^{2n}X_{i_k}\right) 
   = \sum_{P \in {\cal P}} \prod_{r=1}^n E_{\pi}(X_{i_{k_r}}X_{i_{l_r}})
   = \sum_{P \in {\cal P}} \prod_{k=1}^n C_{i_{k_r},i_{l_r}} ~.
\label{eqn:Wick}  \end{equation}
As an example, for $d=1$ and $2n=6$, $X \sim {\cal N}(0,\sigma^2)$, there are $5\cdot 3 = 15$ 
parings, so
\begin{equation}
E(X^6) = 15 \left(\sigma^2\right)^3 = 15 \sigma^6  ~.
\label{eqn:EX^6}  \end{equation}

To apply Wick's formula to our results,
we use the simplified notation 
$F_{ijk} =\partial_{x_i}\partial_{x_j}\partial_{x_k}F(0)$, 
and write
\[
   C_3 = \frac{1}{6}\sum_{ijk} F_{ijk} x_ix_jx_k ~,
\]
and
\[
   C_3^2 = \frac{1}{36}\sum_{ijklmn} F_{ijk}F_{lmn} x_ix_jx_k x_lx_mx_n ~.
\]
There are two kinds of pairings.
One kind pairs one of the indices $\left\{i,j,k\right\}$ with another of the $\left\{i,j,k\right\}$.
This forces one of the $\left\{l,m,n\right\}$ to be paired with another, and the unpaired index from 
$\left\{i,j,k\right\}$ to be paired with the unpaired index from $\left\{l,m,n\right\}$.
An example of this kind of pairing is 
\[
    P = \left\{ \left\{i,k\right\}, \left\{j,n\right\}, \left\{l,m\right\}\right\} ~.
\]
There are nine such pairings, since the unpaired index from each triple is arbitrary.
The expectations are all equal because $F_{ijk}$ is a symmetric function of its indices.
The other kind of pairing has each of the $ \left\{i,j,k\right\}$ paired with one of the 
$\left\{l,m,n\right\}$.
An example of this kind of pairing is 
\[
    P = \left\{ \left\{i,m\right\}, \left\{j,n\right\}, \left\{k,l\right\}\right\} ~.
\]
There are six such pairings, since $i$ is paired with one of the three $ \left\{l,m,m\right\}$, 
then $j$ with one of the remaining two, then $k$ with the last one.
The expectations are again equal.
Altogether
\[
   E_{\pi}(C_3^2) 
     =  \frac{1}{36}\sum_{ijklmn} F_{ijk}F_{lmn} \left( 9 C_{ij}C_{kl}C_{mn}  + 6 C_{il}C_{jm}C_{kn}\right) ~.
\]

This formula simplifies in the special case $H = I$, which implies that $C_{jk} = \delta_{jk}$.
In that case, the $C_{ij}C_{kl}C_{mn}$ terms vanish unless $i=j$, $k=l$ and $m=n$.
The $C_{ij}C_{kl}C_{mn}$ terms give
\[
   \sum_{ikm} F_{iik}F_{mmk} = \left\|\nabla \bigtriangleup F(0) \right\|_{\ell^2}^2  ~.
\]
(For any tensor $A$, we denote the Euclidean 2-norm of all its entries
by $||A||_{\ell^2}$ regardless of the rank of $A~.$) The
$C_{il}C_{jm}C_{kn}$ terms give
\[
   \sum_{ijk} F_{ijk}^2 = \left\| D^3 F(0)\right\|_{\ell^2}^2 ~.
\]
Taken together, these results say that when the Gaussian part of $p$ is
invariant under orthogonal transformations, we have
\[
   E_{\pi}(C_3^2) 
   =  \frac{1}{4} \left\|\nabla \bigtriangleup F(0) \right\|_{l^2}^2 + \frac{1}{6}\left\| D^3 F(0)\right\|_{l^2}^2 ~.
\]
The compact expressions on the right represent the two distinct ways a
quadratic function of the $F_{ijk}$ can be rotationally invariant.

\section{Analysis of random map methods} \label{sec:random-map}
We analyze the simple and symmetrized random map algorithms in the small noise limit $\eps \to 0$.
For the analysis, we use the fact that the random map sampler is affine invariant.
This means that if $M$ is an invertible $d \times d$ matrix and $y=Mx$, 
then the behavior of the random map sampler is identical when applied to $F(x)$ or to $G(x) = F(Mx)$.
Since $H$ is non-degenerate, 
it is possible to choose $M$ so that the Hessian of $G$ is the identity.
Without loss of generality, 
we put $H=I$ in our analysis of random map samplers.
See \cite{G10} for a discussion of the value of affine invariance in practical Monte Carlo.

\subsection{Simple random map}   \label{sec:simple-random-map}

The powers in $Q$ for the simple and symmetrized random map methods come easily.
The simple method has $w = 1+O(\eps^{1/2})$, which the variance lemma (\ref{eqn:Q-var})
turns into $Q=O(\eps)$.
The symmetrized method (\ref{eqn:symmetrized-random-map-weights}) symmetrizes $w$, 
which eliminates the $O(\eps^{1/2})$
term, leaving $w_s = O(\eps)$ and $Q=O(\eps^2)$.
It takes more detailed calculations to find the error constants (\ref{eqn:A-random-map}) 
and (\ref{eqn:A-random-map-symm}).

It is clear that with our assumptions $w$ has an asymptotic expansion in powers of $\eps^{1/2}$
as required by the variance lemma.
For the error constant of the symmetrized method, we need explicit expressions up to 
$O(\eps)$.
We calculate the expansions of the quantities that enter into $w$, then combine them.
We write $a(\xi,\eps) \approx b(\xi,\eps)$ if $a$ and $b$ agree up to order $\eps$.

With our normalization $H=I$, we obtain from (\ref{eqn:F-scaled})
\begin{equation}
    F(x) \approx {\textstyle \frac{1}{2}} |x|^2 + \eps^{1/2}C_3(x) + \eps C_4(x) ~.
\end{equation}
To find an expansion for $\lambda$, we substitute the ansatz
\begin{equation}
  \lambda(\xi) \approx 1 + \eps^{1/2}\lambda_1(\xi) + \eps \lambda_2(\xi) 
  \label{eqn:lambda-ansatz}  
\end{equation}
into (\ref{eqn:stretch-condition}).
We find that
\begin{align}
   {\textstyle \frac{1}{2}} |\xi|^2 
     &\approx {\textstyle \frac{1}{2}}\lambda^2(\xi) |\xi|^2 
         + \eps^{1/2}\lambda(\xi)^3 C_3(\xi) + \eps\lambda(\xi)^4 C_4(\xi), \nonumber\\
     &\approx {\textstyle \frac{1}{2}} |\xi|^2 + \eps^{1/2}\lambda_1 |\xi|^2 
              + \eps \left[ {\textstyle \frac{1}{2}}\lambda_1^2 + \lambda_2 \right] |\xi|^2,   \nonumber\\
       &\;\;\;\;\;  + \eps^{1/2}\left[ 1 + 3\eps^{1/2}\lambda_1\right] C_3(\xi)
          + \eps C_4(\xi) ~.
\end{align}
Collecting terms of $O(\eps^{1/2})$ gives
\[
    0 = \lambda_1(\xi) |\xi|^2 + C_3(\xi) ~,
\]
which can be rearranged to 
\begin{equation}
  \lambda_1(\xi) = \frac{-C_3(\xi)}{|\xi|^2} ~.
\label{eqn:lam1}  \end{equation}
The $O(\eps)$ equation is
\[
   0 = \left[ {\textstyle \frac{1}{2}}\lambda_1(\xi)^2 + \lambda_2(\xi) \right] |\xi|^2 
      + 3 \lambda_1(\xi) C_3(\xi) + C_4(\xi) ~.
\]
Solving for $\lambda_2$ yields
\begin{equation}
    \lambda_2(\xi) = \frac{5}{2}\frac{C_3(\xi)^2}{|\xi|^4} - \frac{C_4(\xi)}{|\xi|^2} ~.
\label{eqn:lam2}  \end{equation}

We now expand the weights (\ref{eqn:w-random-map}).  For the
denominator, we compute the gradient of $F$:
\[
   \nabla F(\xi) \approx \xi + \eps^{1/2}\nabla C_3(\xi) + \eps \nabla C_4(\xi) ~.
\]
Since $C_3(\xi)$ is homogeneous of degree $3$, we have 
$\nabla C_3(\lambda\xi) = \lambda^2 \nabla C_3(\xi)$,
and Euler's identity gives $\xi^t \nabla C_3(\xi) = 3 C_3(\xi)$.
Therefore,
\begin{align}
  \xi^t  \nabla F(\lambda(\xi)\xi) &\approx \lambda(\xi) |\xi|^2 
  + \eps^{1/2} \lambda(\xi)^2 \xi^t \nabla C_3(\xi) 
  + \eps \xi^t \nabla C_4(\xi) \nonumber\\
 &\approx |\xi|^2 + \eps^{1/2}\left[ \lambda_1(\xi) |\xi|^2 + 3 C_3(\xi) \right] \nonumber\\
 & \;\;\;\;\;\;\;\;\;\;\;\;\;
+    \eps \;\left[ \lambda_2(\xi)|\xi|^2 + 6 \lambda_1(\xi) C_3(\xi) + 4 C_4(\xi) \right] \nonumber\\
 &\approx |\xi|^2 + \eps^{1/2}2C_3(\xi)
     + \eps \left[ 3 C_4(\xi) -\,\frac{7}{2} \,\frac{C_3(\xi)^2}{|\xi|^2} \right] ~. \label{eqn:w-bottom}
\end{align}
For the numerator in (\ref{eqn:w-random-map}), use the identity
\[
  (1+\alpha)^{d-1} = 1 + (d-1)\alpha + \frac{1}{2}(d-1)(d-2)\alpha^2 + O(\alpha^3) ~,
\]
to obtain
\begin{align}
   \lambda^{d-1} 
      &\approx \left( 1 + \eps^{1/2}\lambda_1 + \eps\lambda_2 \right)^{d-1} \nonumber\\
      &\approx 1 + \eps^{1/2} (d-1)\lambda_1
          + \eps \left[ (d-1)\lambda_2 + {\textstyle \frac{1}{2}}(d-1)(d-2)\lambda_1^2 \right] \nonumber\\
      &\approx 1 + \eps^{1/2} \frac{(1-d) C_3(\xi)}{|\xi|^2}
         + \eps (d-1)\left[ \frac{d+3}{2} \,\frac{C_2(\xi)^2}{|\xi|^4} 
                           - \frac{C_4(\xi)}{|\xi|^2} \right] ~. \label{eqn:w-top}
\end{align}

We use (\ref{eqn:w-top}) and (\ref{eqn:w-bottom}) to evaluate $w$ to order $\eps$:

\[
  w(\xi) \approx \frac{\displaystyle \left\{ 1 + \eps^{1/2} \frac{(1-d) C_3(\xi)}{|\xi|^2}
         + \eps (d-1)\left[ \frac{d+3}{2} \,\frac{C_2(\xi)^2}{|\xi|^4} 
                           - \frac{C_4(\xi)}{|\xi|^2} \right]\right\} |\xi|^2 }
          {\displaystyle|\xi|^2 + \eps^{1/2}2C_3(\xi)
     + \eps \left[ 3 C_4(\xi) -\,\frac{7}{2} \,\frac{C_3(\xi)^2}{|\xi|^2} \right]} ~.
\]
This has the form 
\begin{equation}
   w \approx  \frac{1+\eps^{1/2}a + \eps b}
         {1+\eps^{1/2}c + \eps d} \approx 1 + \eps^{1/2} (a-c) + \eps \left( c^2-d-ac+b\right),
\label{eqn:w-frac-expansion}  \end{equation}
with coefficients
\begin{align*}
a &= \frac{(1-d) C_3(\xi)}{|\xi|^2}         ~,\\
b &=  (d-1)\left[ \frac{d+3}{2} \,\frac{C_3(\xi)^2}{|\xi|^4} 
                           - \frac{C_4(\xi)}{|\xi|^2} \right]          ~,\\
c &= \frac{2 C_3(\xi)}{|\xi|^2}         ~,\\
d &= \frac{3 C_4(\xi)}{|\xi|^2} -\,\frac{7}{2} \,\frac{C_3(\xi)^2}{|\xi|^4}         ~.
\end{align*}
The term of order $\eps^{1/2}$ in (\ref{eqn:w-frac-expansion}) is
\begin{equation}
   a-c = - (d+1) \frac{C_3^2(\xi)}{|\xi|^2} ~.
\label{eqn:w-random-simple}  \end{equation}

This suffices for the error constant for the simple random map method.
The variance lemma formula (\ref{eqn:Q-var}), together with (\ref{eqn:w-random-simple}) gives
\begin{equation}
  Q \approx \eps (d+1)^2E_{\pi}\left(\frac{C_3(\xi)^2}{|\xi|^4} \right) ~.
\end{equation}
The expected value can be evaluated using the Gaussian integral identity~(\ref{eqn:C4}).
Since $C_3^2(\xi)$ is degree $p=6$, we have
\begin{equation}
   Q \approx \eps \frac{(d+1)^2}{(d+2)(d+4)} 
            E_{\pi}\left( C_3(\xi)^2 \right) ~.
\end{equation}
This is the desired result (\ref{eqn:A-random-map}).

\subsection{Symmetrized random map}   \label{sec:symmetrized-random-map}

The analysis of the symmetrized random map requires the $O(\eps)$ term 
in the $w$ expansion (\ref{eqn:w-frac-expansion}).
The result is
\[
  \frac{(d+2)(d+4)}{2}\frac{C_3(\xi)^2}{|\xi|^4} -
  (d+2)\frac{C_4(\xi)}{|\xi|^2}~.
\]
The $w$ symmetrization formula (\ref{eqn:symmetrized-random-map-weights})
then gives
\[
   w_s(\xi) \approx 1 + 
      \eps \left[ \frac{(d+2)(d+4)}{2}\frac{C_3(\xi)^2}{|\xi|^4} -
  (d+2)\frac{C_4(\xi)}{|\xi|^2} \right]~.
\]
The variance lemma (\ref{eqn:Q-var}) then implies that
\begin{displaymath}
  Q = ~\eps^2~\mbox{var}_{\pi}\Big(\frac{(d+2)(d+4)}{2}\frac{C_3(\xi)^2}{|\xi|^4}
  - (d+2)\frac{C_4(\xi)}{|\xi|^2}\Big)~.
\end{displaymath}
We now expand the above, and rearrange the terms for direct comparison
with the result~(\ref{eqn:A-linear-map-symm}) for the symmetrized linear
map:
\begin{align}
  Q =& ~\eps^2\underbrace{\mbox{var}_{\pi}\Big(\frac{(d+2)(d+4)}{2}\cdot\frac{C_3(\xi)^2}{|\xi|^4}\Big)}_{\mbox{I}}\nonumber\\
  &-\eps^2\underbrace{\mbox{cov}_{\pi}\Big(\frac{(d+2)(d+4)}{2}\cdot\frac{C_3(\xi)^2}{|\xi|^4},
    (d+2)\frac{C_4(\xi)}{|\xi|^2}\Big)}_{\mbox{II}}\nonumber\\
    &+ \eps^2\underbrace{\mbox{var}_{\pi}\Big((d+2)\frac{C_4(\xi)}{|\xi|^2}\Big)}_{\mbox{III}}~.
\end{align}
Consider term I.  We have
\begin{align}
  \mbox{I} &=
  \mbox{var}_{\pi}\Big(\frac{(d+2)(d+4)}{2}\frac{C_3(\xi)^2}{|\xi|^4}\Big) \nonumber\\[2ex]
  &= \frac{(d+2)^2(d+4)^2}{4} E_{\pi}\Big(\frac{C_3(\xi)^4}{|\xi|^8}\Big)
  -\Big[\frac{(d+2)(d+4)}{2} E_{\pi}\Big(\frac{C_3(\xi)^2}{|\xi|^4}\Big)\Big]^2~.
\end{align}
Using (\ref{eqn:C4}) and a direct generalization of it,
we get
\begin{align}
  E_{\pi}\Big(\frac{C_3(\xi)^2}{|\xi|^4}\Big)
    &= \frac{E_{\pi}\Big(C_3(\xi)^2\Big)}{(d+2)(d+4)}, \\
  E_{\pi}\Big(\frac{C_3(\xi)^4}{|\xi|^8}\Big)
    &= \frac{E_{\pi}\Big(C_3(\xi)^4\Big)}{(d+4)(d+6)(d+8)(d+10)} ~.
\end{align}
Thus,
\begin{align}
  \mbox{I}
  &= \frac{(d+2)^2(d+4)^2 E_{\pi}\Big(C_3(\xi)^4\Big)}{4(d+4)(d+6)(d+8)(d+10)}
  -\Big[\frac12 E_{\pi}\Big(C_3(\xi)^2\Big)\Big]^2 \nonumber\\[2ex]
  &= \frac14\mbox{var}_{\pi}\big(C_3(\xi)^2\big)
  - \frac14\cdot\Big\{1-\frac{(d+2)^2(d+4)^2}{(d+4)(d+6)(d+8)(d+10)}\Big\}E_{\pi}\big(C_3(\xi)^4\big)~.
\end{align}
Similarly, we have
\begin{align}
  \mbox{II} &= \frac12\mbox{cov}_{\pi}\Big(C_3(\xi)^2,C_4(\xi)\Big)
  - \frac12\cdot\Big\{1-\frac{(d+2)^2(d+4)}{(d+4)(d+6)(d+8)}\Big\}E_{\pi}\big(C_3(\xi)^2C_4(\xi)\big),\\[2ex]
  \mbox{III} &= \mbox{var}_{\pi}\big(C_4(\xi)\big) -
  \Big\{1-\frac{(d+2)^2}{(d+4)(d+6)}\Big\}E_{\pi}\big(C_4(\xi)^2\big)~.
\end{align}

When $d$ is sufficiently large, we can rewrite the expressions for I-III
in a more concise way:
\begin{align}
  \mbox{I} &= \frac14\mbox{var}_{\pi}\big(C_3(\xi)^2\big)
  - \Big(\frac{4}{d} + O\big(1/d^2\big)\Big)E_{\pi}\big(C_3(\xi)^4\big),\\[1ex]
  \mbox{II} &= \frac12\mbox{cov}_{\pi}\Big(C_3(\xi)^2,C_4(\xi)\Big)
  - \Big(\frac{5}{d} + O\big(1/d^2\big)\Big)E_{\pi}\big(C_3(\xi)^2C_4(\xi)\big),\\[1ex]
  \mbox{III} &= \mbox{var}_{\pi}\big(C_4(\xi)\big) -
  \Big(\frac{6}{d} + O\big(1/d^2\big)\Big)E_{\pi}\big(C_4(\xi)^2\big)~.
\end{align}
This leads to
\begin{equation}
  Q = \eps^2\left[\mbox{var}_{\pi}\Big(\frac12C_3(\xi)^2-C_4(\xi)\Big) +
    c_d\cdot K\right] + O(\eps^3) ~,
\end{equation}
where $c_d = O(1/d)$ and $K$ is a combination of moments of $C_3$ and
$C_4~.$ This verifies the stated result (\ref{eqn:A-random-map-symm}).
We see that for $d\gg1~,$ the variance of the symmetrized random map
method approaches that of the symmetrized linear map.  The above also
shows that in low dimensions, the variance of the symmetrized random map
may be smaller than that of the symmetrized linear map, though exactly
how much depends on the degree of correlation between $C_3(\xi)^2$ and
$C_4(\xi)~.$

The error term of the symmetrized methods in (\ref{eqn:A-linear-map-symm})
and (\ref{eqn:A-random-map-symm}) can in principle also be evaluated using Wick's formula, 
however the calculations are much more involved.
We illustrate how to use Wick's formula for the error terms of
the symmetrized methods with an example.

\section{Computational experiments}   \label{sec:computations}

We present computational experiments that confirm the theoretical error
analysis above for small $\eps$ and suggest what may happen when $\eps$
is not so small.  We use two test problems.  One is a nonlinear random
walk whose dimension is arbitrary.  This allows us to see how the
samplers' performance depends on the dimension.  We see that the
samplers perform worse in higher dimension, but they are still quite
useful in dimensions of practical interest.  In the other example we
apply the algorithms to a data assimilation problem with the 
``Lorenz~'63'' model \cite{L63}.  The goal is is to sample the posterior
distribution of the initial conditions in the presence of noisy
observations of the state at later time.

\subsection{Non-linear random walk}   \label{sec:random-walk}

Consider a non-Gaussian random walk tied at the start and free at the end.
The random variable is $X = (X_1,\ldots,X_N)$, with $X_0\equiv0$ implicitly.
The Gaussian random walk potential is 
\begin{equation}
  x^tHx = \sum_{k=0}^{N-1} \left( x_{k+1} - x_k \right)^2 
  \label{eqn:Gauss-random-walk}
\end{equation}
where $x_0=0$.
We make the walk non-Gaussian by adding cubic and quartic terms to the potential energy.
The nonlinear parts are discretizations of a nonlinear energy functional
\begin{equation}
  C_3(x) = \alpha \sum_{k=0}^{N-1} \left( x_{k+1} - x_k \right)^3 ~,
\end{equation}
and 
\begin{equation}
  C_4(x) = \beta \sum_{k=0}^{N-1} \left( x_{k+1} - x_k \right)^4 ~.
\end{equation}
The coefficients $\alpha$ and $\beta$ would be called ``coupling
constants'' in field theory, and both are set to 1 in our numerical
experiments below.  In the Gaussian measure determined by
(\ref{eqn:Gauss-random-walk}), the increments $\left( X_{k+1} - X_k
\right)$ are independent standard normal random variables.  We can use
this to calculate
\begin{equation}
  E_{\pi}\!\left( C_3(X)^2\right) = 
  \alpha^2\sum_{jk} E_{\pi}\left( \left( X_{j+1} - X_j \right)^3\left( X_{k+1} - X_k \right)^3 \right)~.
  \label{eqn:C-sum}
\end{equation}
The terms on the right hand side with $j \neq k$ vanish because the
increments are independent.  The terms with $j=k$ satisfy, using Wick's
formula (\ref{eqn:EX^6}),
\begin{equation}
  E_{\pi}\!\left( \left( X_{k+1} - X_k \right)^6 \right) = E_{{\cal N}(0,1)}\!\left(Z^6\right) = 15 ~,
\end{equation}
so
\begin{equation}
  E_{\pi}\!\left( C_3(x)^2 \right) = 15 \alpha^2 N ~.
\end{equation}
Thus, the simple linear method for this problem has the quality measure,
 see (\ref{eqn:A-linear-map}),
\begin{equation}
\label{eq:Scaling_Ex1_LM}
  Q \approx \eps 15 \alpha^2 N ~.
\end{equation}
The simple random map quality measure (\ref{eqn:A-random-map}) is slightly less:
\begin{equation}
  Q \approx \eps 15 \alpha^2 \frac{N(N+1)^2}{(N+2)(N+4)} ~.
\label{eqn:rand-map-nonlin-walk}  \end{equation}

It is tedious but straightforward to calculate the error constant for the 
symmetrized methods.
We need  
\[
  \mbox{var}_{\pi}(C_4 - \frac{1}{2}C_3^2) = 
    E_{\pi}\left([C_4 - {\textstyle \frac{1}{2}}C_3^2]^2\right)
  - \left(E_{\pi}(C_4 - {\textstyle \frac{1}{2}}C_3^2)\right)^2 ~.
\]
The first part is
\[
     E_{\pi}\left([C_4 - {\textstyle \frac{1}{2}}C_3^2]^2\right)
  =  E_{\pi}\left(C_4^2\right) 
    -E_{\pi}\left( C_4C_3^2\right)
   + \frac{1}{4}E_{\pi}\left(C_3^4\right)~.
\]
We evaluate these three using Wick identities, first
\begin{align*}
   E_{\pi}\left(C_4^2\right)  &= 
   \beta^2\sum_{jk} E_{{\cal N}(0,1)}\left((X_{j+1}-X_j)^4(X_{k+1}-X_k)^4\right) \\
 &=  \beta^2\Bigl[ \sum_{j\neq k}E\left((X_{j+1}-X_j)^4(X_{k+1}-X_k)^4\right) \\
 &\;\;\;+ \sum_{j= k}E\left((X_{j+1}-X_j)^4(X_{k+1}-X_k)^4\right)\Bigr] \\
 &=  \beta^2\Bigl[(N^2 - N)E\left((X_2-X_1)^4(X_3-X_2)^4\right)\\
 &\;\;\;+ N E\left((X_2-X_1)^8\right)\Bigr] \\
 &=  3 \cdot 3\; \beta^2 N^2 + O(N) \; .
\end{align*}
We write numbers in factored form, as in  $3\cdot 3$ instead of $9$, for 
clarity.
The second term is
\begin{align*}
   E_{\pi}\left( C_4C_3^2\right)  &= 
   \alpha^2 \beta\sum_{jkl} E_{{\cal N}(0,1)}\left((X_{j+1}-X_j)^4(X_{k+1}-X_k)^3(X_{l+1}-X_l)^3\right) \\
 &=  \alpha^2\beta\Bigl[ \sum_{j\neq (k=l)}E\left((X_{j+1}-X_j)^4(X_{k+1}-X_k)^3(X_{l+1}-X_l)^3\right) \\
 &\;\;\;+ \sum_{j= k=l}E\left((X_{j+1}-X_j)^4(X_{k+1}-X_k)^3(X_{l+1}-X_l)^3\right)\Bigr] \\
 &=  \alpha^2\beta\Bigl[(N^2 - N)E\left((X_2-X_1)^4(X_3-X_2)^6\right)\\
 &\;\;\;+ N E\left((X_2-X_1)^{10}\right)\Bigr] \\
 &= 3 \cdot5 \cdot 3\; \alpha^2\beta N^2 + O(N) \; .
\end{align*}
The factor of 3 in the third term is for the three possibilities
$(j=k)\neq(l=m)$, and $(j=l)\neq(k=m)$, and $(j=m)\neq(k=l)$:
\begin{align}
   E_{\pi}\left(C_3^4\right)  &= 
   \alpha^4\sum_{jklm} E\left((X_{j+1}-X_j)^3(X_{k+1}-X_k)^3(X_{l+1}-X_l)^3(X_{m+1}-X_m)^3\right)\nonumber \\
 &=  \alpha^4\Bigl[ \;3\sum_{(j=k)\neq (l=m)}
         E\left((X_{j+1}-X_j)^3(X_{k+1}-X_k)^3(X_{l+1}-X_l)^3(X_{m+1}-X_m)^3\right)\nonumber \\
 &\;\;\;+ \sum_{j= k=l=m}E\left((X_{j+1}-X_j)^3(X_{k+1}-X_k)^3(X_{l+1}-X_l)^3(X_{m+1}-X_m)^3\right)\Bigr] \nonumber\\
 &=  \alpha^4\Bigl[3(N^2 - N)E\left((X_2-X_1)^6(X_3-X_2)^6\right)\nonumber\\
 &\;\;\;+ N E\left((X_2-X_1)^{12}\right)\Bigr] \nonumber\\
 \label{eq:Scaling_Ex1_SLM}
 &= 3 \cdot(5 \cdot 3)^2\;\alpha^4 N^2+O(N) \; .
\end{align}
Adding these gives 
\[
      E_{\pi}\left([C_4 - {\textstyle \frac{1}{2}}C_3^2]^2\right)
   =  N^2 \left( \beta^2\cdot 3^2 - \alpha^2 \beta\cdot 3 \cdot 5 \cdot 3 
             + {\textstyle\frac{1}{4}}\alpha^4 \cdot 3 \cdot(5 \cdot 3)^2\right) + O(N) ~.
\]
A simpler calculation shows that $E_{\pi}(C_4 - {\textstyle \frac{1}{2}}C_3^2) = O(N)$.
Subtracting the terms finally gives
\[
     \mbox{var}_{\pi}(C_4 - \frac{1}{2}C_3^2) 
   = {\textstyle\frac{1}{4}}\alpha^4 N^2\cdot 2 \cdot(5 \cdot 3)^2 + O(N)
   = \frac{225\alpha^4N^2}{2}+ O(N) ~.
\]
It is now clear that the simple methods have 
error coefficients proportional to $\eps N$,
and the symmetrized methods have 
error coefficients proportional to $(\eps N)^2$.

We perform numerical experiments
and vary $N$ and $\eps$.
In these experiments,
we approximate the expected values
in the quality measure $Q$ by
averages over $10^4$ samples.
We protect the computations against
over- and underflow as follows.
Instead of saving the weight, we 
save the logarithm of the weight
of each sample.
This is straightforward for the linear map.
For the symmetrized linear map, 
we use
\begin{align}
  w_{slm}(x)&\propto w_{lm}(x)+w_{lm}(-x) \nonumber\\
  &= w_{lm}(x)\left(1+\frac{w_{lm}(-x)}{w_{lm}(x)}\right),
\end{align}
where $w_{lm}$ is the weight of the simple linear map and $w_{slm}$ that
of the symmetrized linear map.  We then compute
\begin{align}
  \log w_{slm}(x) =& \log(w_{lm}(x)) + \log\left(1+\frac{w_{lm}(-x)}{w_{lm}(x)}\right),\nonumber\\
  =& \log(w_{lm}(x)) + \log\left(1+\exp(F(x)-F(-x))\right).
\end{align}
For the random map we save the log of the weight
\begin{equation}
  \log w_{rm}(x) = (d-1)\log(|\lambda(x)|)+\log(\xi^tH\xi) - \log\left(\xi^t\nabla_x F(\lambda(x)\xi)\right).
\end{equation}
For the symmetrized random map,
the log of the weight is
\begin{align}
  \log w_{srm}(x) &= \log(w_{rm}(x)), \nonumber \\
  &+\log\left(1+\left(\frac{\lambda(-x)}{\lambda(x)}\right)^{d-1}
  \frac{\xi^t\nabla_xF(-\lambda(-x)\xi)}{{\xi^t\nabla_xF(\lambda(x)\xi)} }\right).
\end{align}
Once we have computed the logarithms of the weights for each sample, 
we subtract the maximum value of the logarithms of the weights,
then exponentiate, then normalize.

\begin{figure}[tb]
\includegraphics[width=0.5\textwidth]{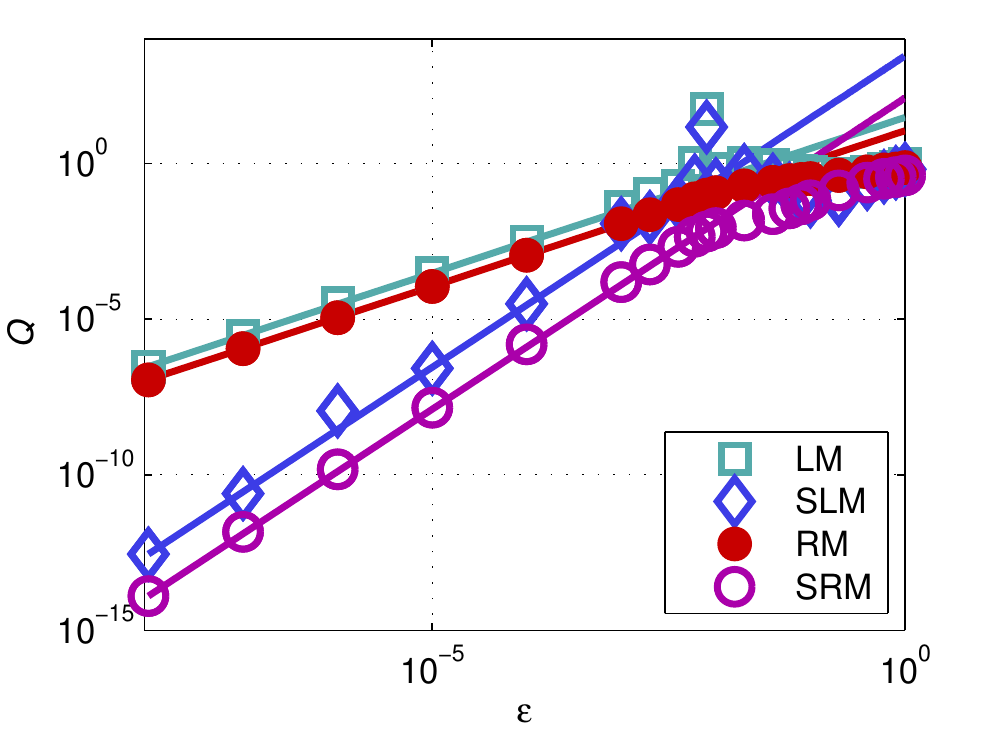}
\includegraphics[width=0.5\textwidth]{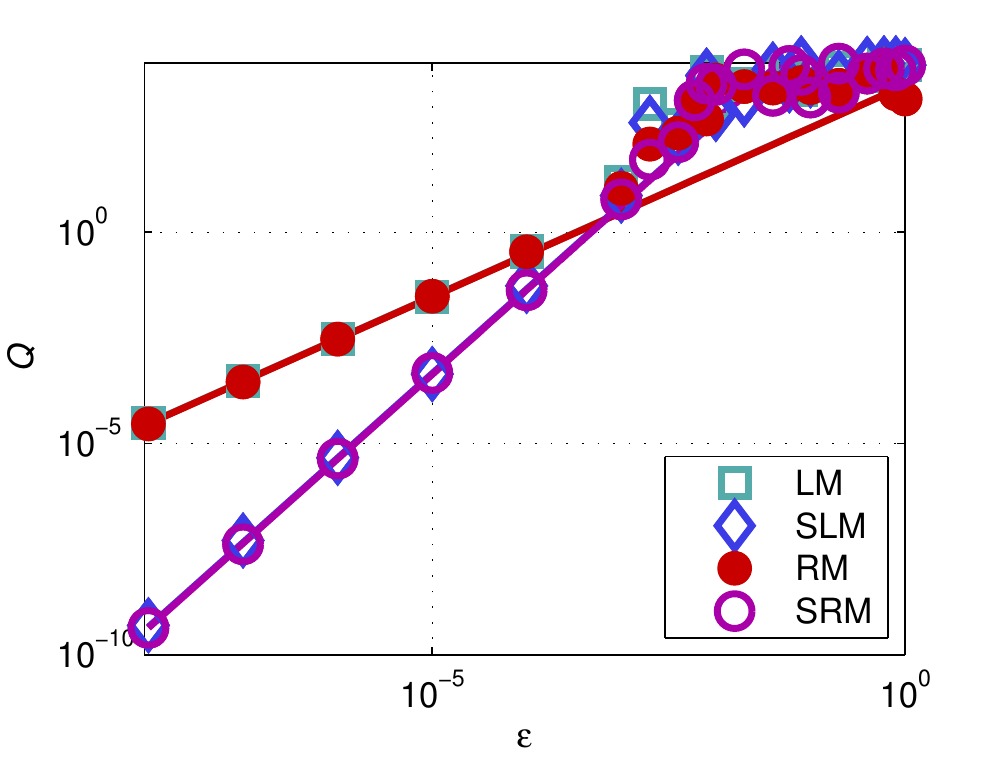}
\caption{Sampling nonlinear random walks.  Left: $N=2$. Right: $N=200$.
Turquoise squares: linear map (LM). Blue diamonds: symmetrized linear map (SLM).
Red dots: random map. Purple circles: symmetrized random map.
The  turquoise and red lines have slope one,
the purple and blue lines have slope two.
The turquoise line is as in (\ref{eq:Scaling_Ex1_LM}),
the red line as in (\ref{eqn:rand-map-nonlin-walk}),
and the purple line on the left is as in (\ref{eq:Scaling_Ex1_SLM}).}
\label{fig:Example1}
\end{figure}

The left panel of Figure~\ref{fig:Example1} shows
$Q$ as a function of $\eps$
for $N=2$, and the right panel for $N=200$.
The dots, circles, squares and diamonds
are values of $Q$ computed from the samples,
the lines have slope one or two,
and are there to illustrate the ``order'' of the method.
Specifically, the turquoise line is as in (\ref{eq:Scaling_Ex1_LM}),
the red line as in (\ref{eqn:rand-map-nonlin-walk}),
and the purple line on the left is as in (\ref{eq:Scaling_Ex1_SLM}).
The numerical results confirm our asymptotic
expansions for sufficiently small $\eps$.
We have made similar observations for 
other values of $N\leq 1000$.
Specifically, for $N=2$, we observe that the numerical results
agree with the predicted values for relatively
large $\eps$ (up to $\eps\approx 0.01$).
For $\eps\geq 10^{-3}$, the
linear map method, the random map method,
and the symmetrized linear map method
are similarly good (as measured by $Q$).
All four methods are doing equally well 
when $\eps$ becomes even larger.
Moreover, all four methods can
be useful in this problem, in the sense
that $Q$ is ``not too large,'' even when
$\eps$ is close to 1.

We observe in the numerical
experiments with $N=200$ that the random map
loses its advantage over the linear map
when $N$ becomes large. 
This is true for the simple and symmetrized
versions of these methods.
We observe that the results of our experiments
agree with our predictions for $\eps\leq 10^{-3}$.
For larger $\eps$, all methods perform
poorly and yield a large $Q\gg 1$
for $\eps\geq 10^{-3}$.

\begin{figure}[tb]
\begin{center}
\includegraphics[width=.6\textwidth]{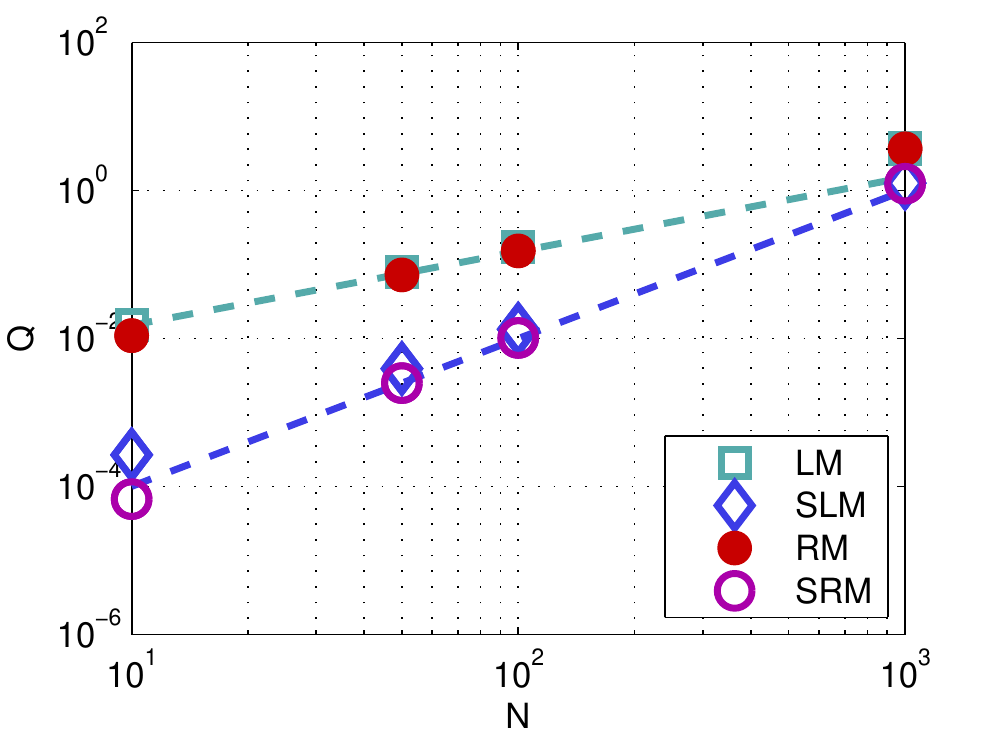}
\end{center}
\caption{Scaling of $Q$ with $N$. 
Turquoise squares: linear map (LM). Blue diamonds: symmetrized linear map (SLM).
Red dots: random map. Purple circles: symmetrized random map.
The  turquoise line has slope one,
the blue line has slope two.}
\label{fig:Example1_Scaling_N}
\end{figure}

Figure \ref{fig:Example1_Scaling_N} illustrates the scaling
of $Q$ with $N$, as computed by Wick's formula.
Shown is $Q$ as a function of $N$ for the various methods.
As  predicted by the theory, we observe that
the symmetrized methods have leading error terms
proportional to $(\eps N)^2$, 
and that the simple, unsymmetrized methods have 
error terms proportional to $\eps N$.

\subsection{Lorenz '63}
We consider estimating the initial conditions of the Lorenz '63 \cite{L63}
equations
\begin{equation}
  \frac{dx}{dt}=\sigma(y-x),\quad  
  \frac{dy}{dt}=x(\rho-z)-y,\quad
  \frac{dz}{dt}=xy-\beta z,
  \label{eqn:lorenz63}
\end{equation}
where $\sigma=10$, $\beta = 8/3$ and $\rho = 28$, from noisy
measurements of $x$, $y$ and $z$ at time $T$:
\begin{equation}
  d = (x(T),y(T),z(T))^t+v.
\end{equation}
The Gaussian random variable $v\sim\mathcal{N}(0,\eps I_3)$ models
measurement noise.  
The above ordinary differential equations (ODE) are solved with the
Matlab routine {\tt ode45}.
The prior for the initial conditions is Gaussian with mean
\begin{equation}
  \mu_0 = (3.6314,\;  6.6136,\;    10.6044)^t,
\end{equation}
and covariance $P_0=\eps I_3$.
The conditional random variable $x_0|d$ thus has the pdf
$p(x_0|d)=\exp(-F(x_0)/\eps)$, where
\begin{equation}
  F(x_0) = \frac{1}{2}\left( (d-h(x_0))^t (d-h(x_0))+  (\mu_0-x_0)^t(\mu_0-x_0)\right),
\end{equation}
so that this problem corresponds to a ``small noise'' situation. 
Here $x_0$ is shorthand notation for the vector $(x(0),y(0),z(0))^t$,
and $h(x_0)$ is the {\tt ode45} solution of the ODEs at time $T$.
The initial conditions we use to generate the synthetic data for our 
numerical experiments is 
\begin{equation*}
	x_{0,\text{true}}= \mu_0 +0.5\,(\sqrt{\eps}, -\sqrt{\eps},\sqrt{\eps})^t
\end{equation*}

We generate samples of $x_0|d$ using the linear and random map methods
described above, and vary $\eps$ and $T$.  The minimization required
by the sampling schemes is done with a quasi-Newton method where all
derivatives are approximated with finite differences.  Similarly, we
approximate the Hessian at the minimum via finite differences.

\begin{figure}[tb]
\begin{center}
\includegraphics[width=.6\textwidth]{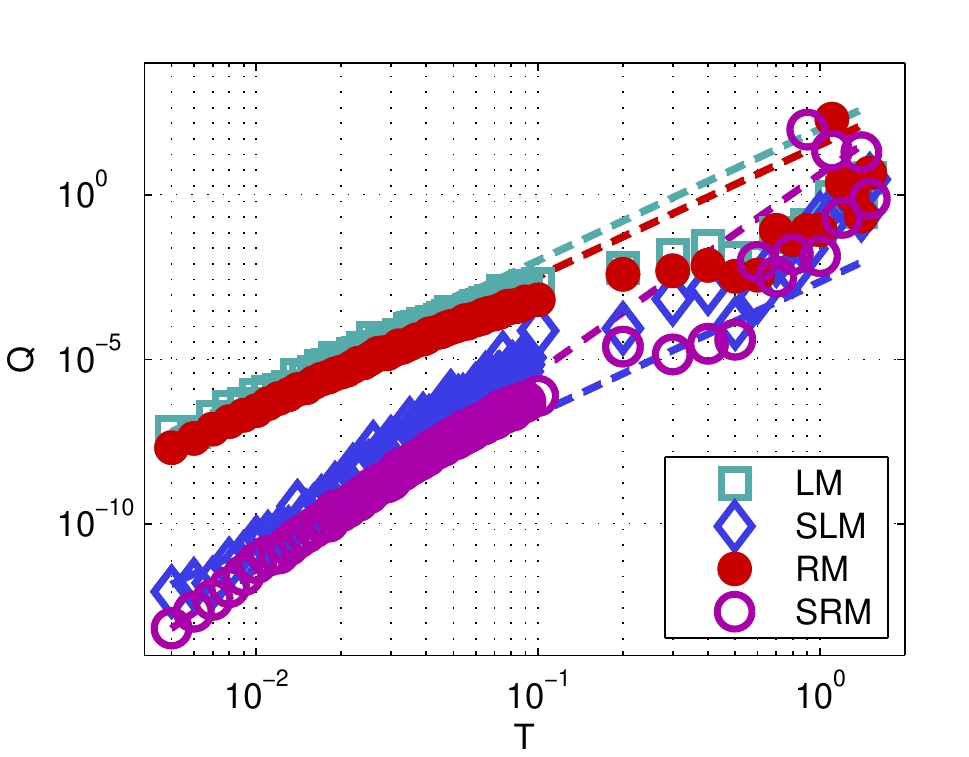}
\caption{Estimating initial conditions of the Lorenz '63 equations.
The parameter $\eps=1$ is constant and the time $T$
at which data are collected is varied.
Turquoise squares: linear map (LM). Blue diamonds: symmetrized linear map (SLM).
Red dots: random map. Purple circles: symmetrized random map.
The  turquoise and red lines have slope four,
the blue and purple lines have slope six.}
\label{fig:L63ConstEps}
\end{center}
\end{figure}

We first fix $\eps=1$ and vary $T$,
i.e.~the time when data are collected.
As $T$ becomes larger, the problem
becomes more and more difficult and 
multiple modes can appear \cite{M94,M99}.
Figure \ref{fig:L63ConstEps} shows $Q$
as a function of $T$.
We observe that the symmetrized methods
perform better than the simple versions,
and give a significantly smaller $Q$-value.
For small $T$, the computed values of $Q$
follow a straight line with slope $4$ for
the random and linear maps, and slope $6$
for the symmetrized methods.
For $T\approx 1$, all four methods perform
similarly well (the symmetrization seems to lose
its advantages) and for $T>1$, the methods perform poorly.
This is perhaps because the pdf we attempt to sample
becomes multi-modal and, therefore, is no longer star-shaped.
However, we made no adjustments to address 
multi-modal target densities.

\begin{figure}[tb]
\begin{center}
\includegraphics[width=.6\textwidth]{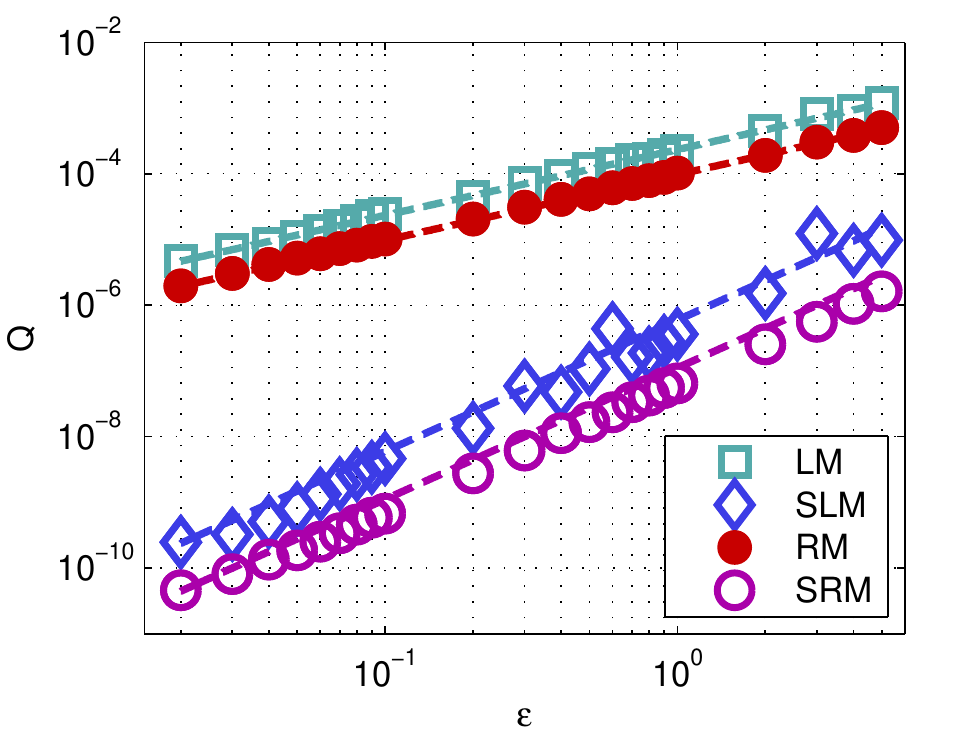}
\caption{Estimating initial conditions of the Lorenz '63 equations.
The data are collected at time $T=0.05$ and the parameter $\eps$
is varied.
Turquoise squares: linear map (LM). Blue diamonds: symmetrized linear map (SLM).
Red dots: random map. Purple circles: symmetrized random map.
The  turquoise and red lines have slope one,
the blue and purple lines have slope two.}
\label{fig:L63_ConstT}
\end{center}
\end{figure}

Next, we fix $T=0.05$ and vary $\eps$.
In this case, the pdf has the functional form
we analyze, and the scenario is analogous
to the ``small noise accurate data'' regime
analyzed in the context of particle filtering 
in \cite{Weare2013}.
Figure \ref{fig:L63_ConstT} shows 
$Q$ as a function of $\eps$.
As in the previous example, 
we find that our numerical experiments
confirm the predicted behavior, 
even if $\eps$ is relatively large.

\section{Conclusion and discussion}\label{sec:Conclusions}

We have performed a small-noise analysis of two implicit sampling
methods, the linear and random map methods.
The analysis shows that the random map method outperforms the linear 
map method in the small noise regime, but this advantage becomes
insignificant in high dimensions.
The simplicity and relative speed of the linear map method thus makes it more
attractive in the limit of small noise.
The analysis further suggests that both methods may be improved by a
symmetrization procedure analogous to antithetic variates.  
We illustrate the theory with numerical examples 
which also suggest that the symmetrized algorithms may
outperform the simple, unsymmetrized algorithms
even when the noise is not so small.

We wish to emphasize two points that are important in practice.  The
first concerns weighted direct samplers as used in particle filtering.
Some methods proposed for practical applications do not have $Q \to
0$, and may even have $Q \to \infty$ as $\eps \to 0$.  For example,
the ``vanilla'' bootstrap particle filter \cite{Gordon93}, which
proposes samples from a proposal distribution that does not take into
account the most recent observation, has $Q \to \infty$
\cite{Weare2012}.  The present samplers all make proposals centered
about the MAP (maximum a-posteriori) point, which takes into account
the most recent observation.  There is much discussion in the
literature of the advantages of doing this
\cite{Chorin2013,Bickel,Doucet,OptimalImportanceFunction,liuchen1995}.

Second, we wish to address the computational cost of the 
algorithms we analyze and propose.
In practice, the cost is roughly proportional to the number
of evaluations of $F$ and its derivatives.
Even our Lorenz `63 example requires an ODE solve to 
evaluate $F$.
The rest of the algorithm is cheap by comparison.

All of our methods start with computing $x_*=\mbox{ argmin }F(x)$.
This requires a  number of evaluations
of $F$ and possibly its derivatives (for numerical optimization).
We also need the Hessian of $F$, either by formulas, 
adjoints, or by finite differences.  
The simple linear map method requires one more evaluation of $F(X)$
per sample. 
The symmetrized linear map method requires two $F$ evaluations.
In particle filter applications, we may want just one sample.
In that case, the optimization is more expensive than sampling.
Other applications may require many samples, in which case the 
cost is roughly the number of samples times the cost for one or
two $F$ evaluations.

The simple random map must solve (\ref{eqn:random-map-stretch}) once for
each sample.  This is one equation in the single unknown, $\lambda$.  It
is normally solved with just a few $F$ evaluations.  As with the linear
map methods, generating one sample using the symmetrized method requires
roughly the work of two samples from the simple version, though by
exploiting the ansatz~(\ref{eqn:lambda-ansatz}) for $\lambda~,$ one can
obtain a good initial guess for finding $\lambda(-\xi)$ based on
$\lambda(\xi)~.$ This may speed up the symmetrized method.

\section*{Acknowledgments}
We gratefully acknowledge the influence of Alexandre Chorin on our work,
and very helpful conversations with him on this material.  Jonathan
Goodman and Kevin Lin thank the Lawrence Berkeley National Laboratory
for facilitating our collaboration on this project.

The work of Jonathan Goodman was supported by the U.S. Department of
Energy, Office of Science, Office of Advanced Scientific Computing
Research, Applied Mathematics Program under contract number
DE-AC02005CH11231 under a subcontract from Lawrence Berkeley National
Laboratory to New York University.  The work of Kevin Lin was supported
in part by the National Science Foundation under grants DMS-1217065 and
DMS-1418775.  The work of Matthias Morzfeld was supported by the
U.S. Department of Energy, Office of Science, Office of Advanced
Scientific Computing Research, Applied Mathematics Program under
contract number DE-AC02005CH11231, and by the National Science
Foundation under grant DMS-1217065.


\bibliographystyle{siam}
\bibliography{implicit}

\begin{thebibliography}{10}

\bibitem{An}
{\sc S.~An and F.~Schorfheide}, {\em Bayesian analysis of {DSGE} models},
  Econometric Reviews, 26 (2007), pp.~113--172.

\bibitem{GordonReview}
{\sc M.~Arulampalam, S.~Maskell, N.~Gordon, and T.~Clapp}, {\em A tutorial on
  particle filters for online nonlinear/non-{G}aussian {B}ayesian tracking},
  IEEE Transactions on Signal Processing, 50 (2002), pp.~174 --188.

\bibitem{atkins}
{\sc E.~Atkins, M.~Morzfeld, and A.~Chorin}, {\em Implicit particle methods and
  their connection with variational data assimilation}, Monthly Weather Review,
  141 (2013), pp.~1786--1803.

\bibitem{Bergmman99}
{\sc N.~Bergman}, ed., {\em Recursive Bayesian estimation: Navigation and
  tracking applications}, Ph.D Dissertation, Linkoping University, Linkoping,
  Sweden, 1999.

\bibitem{Bocquet2010}
{\sc M.~Bocquet, C.~Pires, and L.~Wu}, {\em Beyond {G}aussian statistical
  modeling in geophysical data assimilation}, Monthly Weather Review, 138
  (2010), pp.~2997--3023.

\bibitem{ChorinHald}
{\sc A.~Chorin and O.~Hald}, {\em Stochastic Tools in Mathematics and Science},
  Springer, third~ed., 2013.

\bibitem{Chorin2013}
{\sc A.~Chorin and M.~Morzfeld}, {\em Conditions for successful data
  assimilation}, Journal of Geophysical Research, 118 (2013),
  pp.~11,522--11,533.

\bibitem{Chorin2010}
{\sc A.~Chorin, M.~Morzfeld, and X.~Tu}, {\em Implicit particle filters for
  data assimilation}, Communications in Applied Mathematics and Computational
  Science, 5 (2010), pp.~221--240.

\bibitem{chorintupnas}
{\sc A.~Chorin and X.~Tu}, {\em Implicit sampling for particle filters},
  Proceedings of the National Academy of Sciences, 106 (2009),
  pp.~17249--17254.

\bibitem{Doucet2001}
{\sc A.~Doucet, N.~de~Freitas, and N.~Gordon}, eds., {\em Sequential {M}onte
  {C}arlo {M}ethods in Practice}, Springer, 2001.

\bibitem{Doucet}
{\sc A.~Doucet, S.~Godsill, and C.~Andrieu}, {\em On sequential {M}onte {C}arlo
  sampling methods for {B}ayesian filtering}, Statistics and Computing, 10
  (2000), pp.~197--208.

\bibitem{Fournier2010}
{\sc A.~Fournier, G.~Hulot, D.~Jault, W.~Kuang, W.~Tangborn, N.~Gillet,
  E.~Canet, J.~Aubert, and F.~Lhuillier}, {\em And introduction to data
  assimilation and}, Space Science Review, 137 (2009), pp.~247--291.

\bibitem{G10}
{\sc J.~Goodman and J.~Weare}, {\em Ensemble samplers with affine invariance},
  Communications in Applied Mathematics and Computational Science, 5 (2010),
  pp.~65--80.

\bibitem{Gordon93}
{\sc N.~Gordon, D.~Salmond, and A.~Smith}, {\em Novel approach to
  nonlinear/non-{G}aussian {B}ayesian state estimation}, Radar and Signal
  Processing, IEEE Proceedings F, 140 (1993), pp.~107--113.

\bibitem{HammersleyHandscomb}
{\sc J.~Hammersley and D.~Handscomb}, {\em Monte {C}arlo Methods}, Chapman \&
  Hall, 1964.

\bibitem{KalosWhitlock}
{\sc M.~Kalos and P.~Whitlock}, {\em Monte {C}arlo Methods}, vol.~1, John Wiley
  \& Sons, 1~ed., 1986.

\bibitem{Koopmans}
{\sc L.~Koopmans}, {\em The spectral analysis of time series}, Academic Press,
  1974.

\bibitem{liuchen1995}
{\sc J.~Liu and R.~Chen}, {\em Blind deconvolution via sequential imputations},
  Journal of the American Statistical Association, 90 (1995), pp.~567--576.

\bibitem{LiuChen98}
\leavevmode\vrule height 2pt depth -1.6pt width 23pt, {\em Sequential {M}note
  {C}arlo methods for dynamical systems}, Journal of the American Statistical
  Association, 93 (1998), pp.~1032--1044.

\bibitem{L63}
{\sc E.~Lorenz}, {\em Deterministic nonperiodic flow}, Journal of the
  Atmospheric Sciences, 20 (1963), pp.~130--141.

\bibitem{M99}
{\sc R.~Miller, J.~Carter, and S.~Blue}, {\em Data assimilation into nonlinear
  stochastic models}, Tellus, 51 (1999), pp.~167--194.

\bibitem{M94}
{\sc R.~Miller, M.~Ghil, and F.~Gauthiez}, {\em Advanced data assimilation in
  strongly nonlinear dynamical systems}, Journal of Atmospheric Science, 51
  (1994), pp.~1037--1056.

\bibitem{Morzfeld2011}
{\sc M.~Morzfeld, X.~Tu, E.~Atkins, and A.~Chorin}, {\em A random map
  implementation of implicit filters}, Journal of Computational Physics, 231
  (2012), pp.~2049--2066.

\bibitem{Murray}
{\sc J.~Murray}, {\em Asymptotic Analysis}, Springer Verlag, 1992.

\bibitem{Bickel}
{\sc C.~Snyder, T.~Bengtsson, P.~Bickel, and J.~Anderson}, {\em Obstacles to
  high-dimensional particle filtering}, Monthly Weather Review, 136 (2008),
  pp.~4629--4640.

\bibitem{PeterJan2009}
{\sc P.~van Leeuwen}, {\em Particle filtering in geophysical systems}, Monthly
  Weather Review, 137 (2009), pp.~4089--4114.

\bibitem{Weare2012}
{\sc E.~Vanden-Eijnden and J.~Weare}, {\em Rare event simulation and small
  noise diffusions}, Communications on Pure and Applied Mathematics, 65 (2012),
  pp.~1770--1803.

\bibitem{Weare2013}
\leavevmode\vrule height 2pt depth -1.6pt width 23pt, {\em Data assimilation in
  the low noise, accurate observation regime with application to the kuroshio
  current}, Monthly Weather Review,  (2013).

\bibitem{OptimalImportanceFunction}
{\sc V.~Zaritskii and L.~Shimelevich}, {\em Monte {C}arlo technique in problems
  of optimal data processing}, Automation and Remote Control, 12 (1975), pp.~95
  -- 103.

\end{thebibliography}

\end{document}